\documentclass[11pt,a4paper]{article}
\usepackage[numbers,sort&compress]{natbib}
\usepackage[colorlinks,
            linkcolor=blue,
            anchorcolor=green,
            citecolor=magenta
            ]{hyperref}
\usepackage{mathrsfs,amssymb,amsfonts}
\usepackage{graphicx}
\usepackage{amsmath,amssymb,latexsym,color}
\usepackage[mathscr]{eucal}
\usepackage[numbers]{natbib}
\usepackage{bm}
\usepackage[top=1in, bottom=1in, left=1.25in, right=1.25in]{geometry}
\usepackage{tikz}
\usepackage{verbatim}
\usetikzlibrary{arrows,decorations.pathmorphing,backgrounds,positioning,fit,petri}
\usepackage{setspace}
\makeatletter

\newcommand{\Rmnum}[1]{\expandafter\@slowromancap\romannumeral #1@}
\makeatother
\newtheorem{theorem}{Theorem}[section]

\newtheorem{lemma}[theorem]{Lemma}
\newtheorem{remark}[theorem]{Remark}

\newtheorem{example}[theorem]{Example}

\newcommand{\qed}{\hfill\Box\medskip}
\usepackage{geometry}
\begin{document}
\begin{spacing}{1.5}


\title{\bf  On the K\H{o}nig-Egerv\'{a}ry theorem for $k$-trees \footnote{This research is supported by the National Natural Science Foundation of China(12301453), the Project of Scientific Research Fund of Hunan Provincial Science and Technology Department (2018WK4006)}
}

\author{Pingshan Li\footnote{ E-mail address: lips@xtu.edu.cn}
\\
{\footnotesize \em School of Mathematics and Computational Science, Xiangtan University,
Xiangtan, Hunan 411105, PR China}\\
{\footnotesize Key Laboratory of Intelligent Computing \& Information Processing of Ministry of Education}\\
{\footnotesize Key Laboratory for Computation and Simulation in Science and Engineering}\\
{\footnotesize National Center for Applied Mathematics in Hunan}\\
}

 \date{}
 \date{}
 \maketitle

\begin{abstract}
The famous K\H{o}nig-Egerv\'{a}ry theorem is equivalent to the statement that a graph is bipartite if and only if the matching number equals the vertex cover number for every induced subgraph. Bessy, Ochem, and Rautenbach (J. Graph Theory 91 (2019) 73-87) generalized this theorem to $k$-paths by introducing the class $\mathcal{G}_k$ of graphs where, for every induced subgraph, the maximum number of disjoint paths of order $k$ equals the minimum number of vertices intersecting all such paths. They provided complete structural characterizations of graphs in $\mathcal{G}_k$ for $k \in \{3, 4\}$, and for odd $k$ when no cycle of order less than $k$ is present. Inspired by this nice work, we extend the K\H{o}nig-Egerv\'{a}ry theorem to $k$-trees (i.e., trees with exactly $k$ vertices). We define the class $\mathcal{J}_k$ of graphs where, for every induced subgraph, the maximum number of disjoint trees of order $k$ equals the minimum number of vertices intersecting all such trees, and provide complete structural descriptions of graphs in $\mathcal{J}_k$ for all $k \ge 3$.

\medskip

\noindent {\em Key words:}  Dissociation sets; $k$-path cover, K\H{o}nig-Egerv\'{a}ry theorem, matching.

\medskip
\end{abstract}

\section{Introduction}

An independent edge set in a graph is referred to as a \emph{matching}. The dual concept of a matching is a vertex cover: a vertex subset $U \subseteq V(G)$ such that every edge in $G$ is incident to at least one vertex in $U$. For any graph $G$, let $\nu(G)$ denote its matching number (the size of a maximum matching) and $\tau(G)$ its vertex cover number (the size of a minimum vertex cover). A fundamental observation in graph theory is that $\nu(G) \leq \tau(G)$ holds for all graphs.

The famous K\H{o}nig-Egerv\'ary theorem \cite{egervary,donig} refines this inequality for bipartite graphs, stating that $\nu(G) = \tau(G)$ when $G$ is bipartite. Equivalently, this theorem characterizes bipartite graphs as those where every induced subgraph $H$ satisfies $\nu(H) = \tau(H)$. This equivalence follows from the fact that odd cycles $C_{2k+1}$ violate the equality: for such cycles, $\nu(C_{2k+1}) = k < \tau(C_{2k+1}) = k+1$.

Matching theory can be naturally extended beyond edges to more general subgraph structures through the well-studied packing and covering problems. Given a specified graph class $\mathcal{H}$, the packing problem seeks the maximum number of vertex-disjoint subgraphs of $G$ that are isomorphic to members of $\mathcal{H}$. The dual covering problem aims to find the smallest vertex subset $U \subseteq V(G)$ such that every subgraph of $G$ isomorphic to a member of $\mathcal{H}$ intersects $U$.

Note that the classic matching and vertex cover correspond to the special case where $\mathcal{H} = \{K_2\}$, the graph consisting of a single edge. A common generalization focuses on $\mathcal{H} = \{P_k\}$, where $P_k$ is a path with exactly $k$ vertices. In this setting: A \emph{$k$-matching} is a collection of vertex-disjoint $k$-paths in $G$, with the maximum size denoted $\nu'_k(G)$.
A \emph{$k$-vertex cover} is a vertex set that intersects every $k$-path in $G$, with the minimum size denoted $\tau'_k(G)$.

Between the two parameters $\nu'_k(G)$ and $\tau'_k(G)$, only the latter seems to have received considerable attention in the literature \cite{B2011, B2013, Li2016}. For all integers $k \ge 3$, the computational hardness of the $k$-vertex cover problem was established by Bre\v{s}ar et al.\ \cite{B2011}, and the 3-vertex cover number has been studied extensively \cite{C6, C7, C10, C11,C13, C12}. Notably, a vertex set $X$ of a graph $G$ forms a 3-vertex cover if and only if its complement $V(G)\setminus X$ is a so-called dissociation set, a set of vertices inducing a subgraph with maximum degree at most 1. Related problems involving dissociation sets have recently gained significant attention, as demonstrated in \cite{bock1, bock2, sun, Tu2021}.

Notably, Tu, Zhang and Shi \cite{Tu2021} showed that for every forest $F$ and every integer $k\ge2$, $\nu'_k(F) = \tau'_k(F)$. Bessy, Ochem and Rautenbach \cite{bessy} generalized the K\H{o}nig-Egerv\'{a}ry theorem to long paths by introducing the graph class $\mathcal{G}_k$, consisting of all graphs $G$ such that $\nu'_k(H) = \tau'_k(H)$ for every induced subgraph $H$ of $G$. They provided complete structural characterizations of graphs in $\mathcal{G}_k$ for $k \in \{3, 4\}$, and for odd $k$ when no cycle of order less than $k$ is present.


Since the path $P_k$ is a tree on $k$ vertices, the setting can be further generalized to cases where $\mathcal{H}$ contains all $k$-trees (i.e., trees with exactly $k$ vertices for a fixed integer $k$). Correspondingly, classic graph-theoretic notions of matching and vertex cover extend naturally to their $k$-tree analogs: the generalized $k$-matching and generalized $k$-vertex cover.

\begin{itemize}
  \item A collection of vertex-disjoint $k$-trees in $G$ is referred to as a \emph{generalized $k$-matching}, with the size of a maximum such collection denoted the \emph{generalized $k$-matching number} $\nu_k(G)$.
  \item A \emph{generalized $k$-vertex cover} is a vertex set that intersects every $k$-tree in $G$, with the minimum size of such a set termed the \emph{generalized $k$-vertex cover number} $\tau_k(G)$.
\end{itemize}
It follows directly that \( \nu_k(G) \leq \tau_k(G) \) holds for any graph \( G \). Notably, when \( k = 3 \), the $3$-vertex cover coincides with the generalized $3$-vertex cover, and the $3$-matching coincides with the generalized $3$-matching.

Let $\mathcal{J}_k$ denote the class of all graphs $G$ such that $\nu_k(H)=\tau_k(H)$ holds for every induced subgraph $H$ of $G$. Then, a natural question arises: how can the K\H{o}nig-Egerv\'ary theorem be extended to $k$-trees? In other words, how do we characterize the graphs in $\mathcal{J}_k$?

In 2024, Li and Xu \cite{lps} established that if $F$ is a forest, then $\nu_k(F)=\tau_k(F)$. This result can be restated as follows:

\begin{theorem}{\rm \cite{lps}}\label{tree} Every forest belongs to $\mathcal{J}_k$ for every integer $k\ge2$.
\end{theorem}

When $k=3$, we have $\mathcal{J}_3=\mathcal{G}_3$, so the characterization in \cite{bessy} applies. In this paper, we provide complete structural characterizations of graphs in $\mathcal{J}_k$ for all integers $k\ge3$.

\section{Preliminaries}

Throughout this paper, graphs are finite, simple and undirected; the auxiliary multigraphs may contain loops and parallel edges. In a multigraph, a loop and a pair of parallel edges count as cycles of lengths one and two, respectively. A $k$-tree means a tree of order $k$, not necessarily an induced subgraph. Let $V(G)$ and $E(G)$ represent the vertex set and edge set of a graph $G$, respectively. For any graph $G$ or vertex set $X$, let $n(\cdot)$ denote the cardinality of the underlying vertex set: specifically, $n(G) = |V(G)|$ and $n(X) = |X|$.
A connected graph of order at least $k$ contains a tree of order $k$. Consequently, a vertex set $X$ is a generalized $k$-vertex cover if and only if every component of $G-X$ has fewer than $k$ vertices. We use this observation throughout.
To provide a comprehensive structural characterization of graphs in $\mathcal{J}_k$ for $k \geq 3$, we first introduce a new concept: the $k$-tree insertion operation.

Let $G$ be a graph and let $uv \in E(G)$ be an edge of $G$. The \emph{$k$-tree insertion operation} on edge $uv$ is defined through the following sequential steps:

\begin{itemize}
    \item Remove the edge $uv$ from $G$. 
    \item Introduce a new connected graph $G_{uv}$ with $n(G_{uv}) = k-1$ that is disjoint from $G$.
    \item Select two vertices $x, y \in V(G_{uv})$ (allowing $x = y$ when $u\ne v$; for a loop, require $x\ne y$).
    \item Connect $G_{uv}$ to $G - uv$ by adding edges $ux$ and $vy$.
\end{itemize}

For an integer $k\ge3$, the $k$-tree insertion operation applied to a multigraph $H$, denoted $\text{Tio}_k(H)$, is obtained by performing the $k$-tree insertion operation on every edge of $H$, using mutually disjoint new graphs. For parallel edges, each subscript $uv$ refers to the particular edge under consideration. Specifically, when each added graph $G_{uv}$ is a path and the chosen vertices $x, y$ are the two endpoints of this path, the resulting graph is called the $k$-subdivision of $H$, denoted $\text{Sub}_k(H)$.

The following lemma demonstrates that constructing a graph satisfying $\nu_k(G)=\tau_k(G)$ can be achieved through straightforward means, a result analogous to the analysis in \cite{bessy} for graphs containing a block $B = \text{Sub}_k(H)$ for some multigraph $H$ when $\mathcal{H} = \{P_k\}$. Specifically, if $G = \text{Tio}_k(H)$ where $H$ is a connected multigraph containing at least one cycle, then $\nu_k(G)=\tau_k(G)=n(H)$ for all $k\ge3$. This equality for $G$ alone does not assert that $G\in\mathcal{J}_k$.

To formalize the proof of this lemma, we first introduce a standard graph-theoretic concept: A rooted tree is formed by designating a distinguished vertex in a tree as the root. For a rooted tree $T$, the level of a vertex $u$, denoted by $\ell(u)$, is defined as the distance between $u$ and the root of $T$.

\begin{lemma}\label{multigraph}
Let $k\ge3$. Suppose the graph $G$ contains an induced subgraph $B$ such that $B={{\rm Tio}}_k(H)$ for some connected multigraph $H$ that contains a cycle. Suppose also that every component $K$ of $G-V(H)$ meeting $V(B)\setminus V(H)$ satisfies $\nu_k(K)=0$. Then \begin{description}
       \item[(i)] $\nu_k(G)=\nu_k(G-V(H))+n(H)$;
       \item[(ii)] $\tau_k(G)=\tau_k(G-V(H))+n(H)$.
     \end{description}
\end{lemma}

\begin{proof} (i) Since $H$ is connected and contains a cycle, there exists an edge $e$ incident to some vertex $r$ such that $H-e$ contains a spanning tree $T$ of $H$. Root $T$ at $r$. For each $uv \in E(T)$ with $\ell(u) < \ell(v)$, let $T_{uv}$ be a spanning tree of $B[V(G_{uv}) \cup \{v\}]$. Similarly, let $T_e$ be a spanning tree of $B[V(G_e) \cup \{r\}]$. Define:
\[
\mathcal{A}_0 = \{ T_{uv} \mid uv \in E(T), \ell(u) < \ell(v) \} \cup \{ T_e \}.
\]
The collection $\mathcal{A}_0$ consists of $n(H)$ vertex-disjoint $k$-trees. By the additional hypothesis, every component $K$ of $G-V(H)$ that contains a vertex from $V(B)\setminus V(H)$ satisfies $\nu_k(K)=0$. Therefore, adding $\mathcal{A}_0$ to any generalized $k$-matching of $G-V(H)$ yields a generalized $k$-matching of $G$. It follows that $\nu_k(G)\ge \nu_k(G-V(H))+n(H)$.

Conversely, if $\mathcal{M}$ is a generalized $k$-matching in $G$, then since every $k$-tree in $G$ intersecting $V(B)$ contains a vertex of $H$, the set $\mathcal{M}$ contains at most $n(H)$ $k$-trees intersecting $V(B)$. Removing all such $k$-trees from $\mathcal{M}$ yields a generalized $k$-matching of $G-V(H)$, which implies that
\(
\nu_k(G) \leq \nu_k(G-V(H)) + n(H).
\)

By combining the upper and lower bounds established above, we conclude that
\(
\nu_k(G)= \nu_k(G-V(H)) + n(H).
\)

(ii) If \(X\) is a generalized \(k\)-vertex cover in \(G-V(H)\), then \(X \cup V(H)\) forms a generalized \(k\)-vertex cover in \(G\). This directly yields the inequality
\(
\tau_k(G) \leq \tau_k(G-V(H)) + n(H).
\)

For the converse direction, let \(X\) denote a generalized \(k\)-vertex cover of \(G\). Set \(X'=X\cap V(B)\). Since $X'$ meets every tree in $\mathcal{A}_0$, we have $|X'|\ge n(H)$. By the component hypothesis, every $k$-tree in $G-V(H)$ is disjoint from $V(B)$. Thus $X\setminus X'$ is a generalized $k$-vertex cover of $G-V(H)$, and
\[
\tau_k(G)\ge \tau_k(G-V(H))+n(H).
\]

By combining the upper and lower bounds established above, we conclude that \(\tau_k(G) = \tau_k(G-V(H)) + n(H)\).    $\qed$
\end{proof}

We next establish a straightforward property of the graph family $\mathcal{J}_k$, which is formalized in the lemma below:

\begin{lemma}\label{cycle-lemma}Let $k\ge3$ and let $p$ be a positive integer.
Then \begin{description}
\item[(i)] Every cycle of length at least $k$ in every graph in $\mathcal{J}_k$ has length congruent to $0$ modulo $k$.
\item[(ii)]Let $C_{pk} = u_1u_2\cdots u_{pk}u_1$ be a cycle of length $pk$. A vertex subset $X \subseteq V(C_{pk})$ is a minimum generalized $k$-vertex cover of $C_{pk}$ if and only if $X = \left\{ u_{i+(j-1)k} : j \in [p] \right\}$ for some $i \in [k]$.
\item[(iii)] Suppose $G \in \mathcal{J}_k$, $C$ is a cycle in $G$ with $n(C) \geq k$, and $u, v$ are distinct vertices of $C$ each having at least one neighbor outside $V(C)$. Then $\mathrm{dist}_C(u, v) \equiv 0 \pmod{k}$.
\end{description}
\end{lemma}
\begin{proof} Let $C_n$ be a cycle of length $n\ge k$ in $G\in\mathcal{J}_k$, and let $F=G[V(C_n)]$. Then
\[
\left\lfloor\frac{n}{k}\right\rfloor=\nu_k(C_n)\le\nu_k(F)\le\left\lfloor\frac{n}{k}\right\rfloor.
\]
Every generalized $k$-vertex cover of a cycle must meet every set of $k$ consecutive vertices, and selecting vertices with gaps of at most $k-1$ gives a cover of cardinality $\lceil n/k\rceil$. Thus
\[
\nu_k(F)=\left\lfloor\frac{n}{k}\right\rfloor
\le\left\lceil\frac{n}{k}\right\rceil=\tau_k(C_n)\le\tau_k(F).
\]
Since $F$ is an induced subgraph of $G\in\mathcal{J}_k$, we have $\nu_k(F)=\tau_k(F)$, from which part (i) follows.

Statement (ii) follows directly from two key observations: first, the generalized $k$-vertex cover number of the cycle graph $C_{pk}$ satisfies $\tau_k(C_{pk}) = p$; second, every generalized $k$-vertex cover of $C_{pk}$ must contain at least one vertex from every set of $k$ consecutive vertices on the cycle.

Suppose $G$, cycle $C$, and vertices $u, v$ satisfy the hypotheses of part (iii). By (i), we have $n(C)\equiv 0\mod k$. Let $u' \in V(G) \setminus V(C)$ be a neighbor of $u$, $v' \in V(G) \setminus V(C)$ be a neighbor of $v$, and define $G'$ as the induced subgraph of $G$ with vertex set $V(C) \cup \{u', v'\}$. Then the generalized $k$-matching number of $G'$ satisfies:
\[
\nu_k(G') = \left\lfloor \frac{n(C) + |\{u', v'\}|}{k} \right\rfloor = \frac{n(C)}{k}.
\]
Since $G \in \mathcal{J}_k$ and $G'$ is an induced subgraph of $G$, we have $\tau_k(G') = \nu_k(G')$. Thus,
 \begin{center}
 $\tau_k(C) = \nu_k(C) = \frac{n(C)}{k}=\nu_k(G')=\tau_k(G')$.
\end{center}
It follows that every minimum generalized $k$-vertex cover of $G'$ is also a minimum generalized $k$-vertex cover $X$ of $C$. By part (ii), such a set $X$ must take the form
 \begin{center}
 $X = \left\{ u_{i+(j-1)k} : j \in [p] \right\} = \left\{ x_1, x_2, \dots, x_p \right\}$ for some $i \in [k]$.
\end{center}
Note that removing $X$ from $C$ yields $p$ disjoint paths each of order $k-1$, and there is no $k$-tree in $G' - X$. If $u\notin X$, its component in $C-X$, together with $u'$, is a connected subgraph of order $k$, a contradiction. Thus $u\in X$, and similarly $v\in X$. Thus, $\mathrm{dist}_C(u, v) \equiv 0 \pmod{k}$.
$\qed$
\end{proof}

Finally, we introduce a classic result on 2-connected graphs: the Ear Decomposition Theorem. This theorem establishes a fundamental characterization of 2-connected graphs through a constructive decomposition, which we now define formally.

An open ear of a graph is a path whose distinct endpoints belong to the graph and whose internal vertices do not. An ear decomposition of $G$ is a sequence of edge-disjoint subgraphs $P_0,P_1,\dots,P_s$ whose union is $G$, such that $P_0$ is a cycle and, for each $i\ge1$, $P_i$ is an open ear of $P_0\cup\cdots\cup P_{i-1}$. Ears of length one are allowed. This structural decomposition provides a powerful framework for characterizing 2-connected graphs, as formalized by the following fundamental theorem:

\begin{theorem}{\rm \cite{whitney}}\label{ear}
A graph with at least three vertices is $2$-connected if and only if it has an ear decomposition. Furthermore, every cycle in a $2$-connected graph can serve as the initial cycle in some ear decomposition of the graph.
\end{theorem}

\section{The graphs in $\mathcal{J}_k$}

In this section, we provide a comprehensive structural characterization of graphs in the family $\mathcal{J}_k$ for $k \geq 3$. Specifically, we explicitly characterize their block structures and formalize constraints governing their cut-vertices. We denote the cut-vertices of $G$ contained in $B$ by $x_1,\dots,x_m$. Let $C_B$ be the component of $G$ containing $B$, and let $X_i$ be the union of the vertex sets of all components of $C_B-x_i$ disjoint from $B$. Thus unrelated components of $G$ are not counted in $X_i$. The sets $X_1,\dots,X_m$ are pairwise disjoint. An endblock of $G$ is defined as a block containing at most one cut-vertex of $G$. A block is trivial if it is isomorphic to either $K_1$ or $K_2$.

To formalize our characterization, we first define a new family of graphs that captures these structural properties precisely. Let $\mathcal{H}_k$ denote the family of all graphs $G$ such that every nontrivial block $B$ of $G$ satisfies one of the following conditions:

\begin{description}
  \item[(i)] $B$ is an endblock on exactly $k$ vertices;

  \item[(ii)] $B$ is a block on $t$ vertices where $t < k$, containing $m$ cut-vertices of $G$ with $m \leq k - t + 1$. Moreover, when $m \geq 2$, there exists a choice of $m-1$ sets from $X_1,X_2,\dots,X_m$ whose union has fewer than $k-t+1$ vertices;

  \item[(iii)] $B = \text{Sub}_k(H)$ for some multigraph $H$, and every cut-vertex of $G$ contained in $B$ is a vertex of $H$.
\end{description}

For $m\ge2$, the additional requirement in (ii) is equivalently
\begin{equation}\label{small-condition}
\text{there exists }j\in[m]\text{ such that }\sum_{i\ne j}n(X_i)\le k-t.
\end{equation}
If $B$ satisfies (i) or (ii), there is a vertex $r\in V(B)$ such that
\begin{equation}\label{small-region}
U=V(B)\cup\bigcup_{x_i\ne r}X_i
\quad\text{satisfies}\quad n(U)\le k.
\end{equation}
Here $r$ is the only cut-vertex when there is exactly one, and may be chosen arbitrarily when there is none. The graph $G[U]$ is connected, and no edge joins $U\setminus\{r\}$ to $V(G)\setminus U$. Conversely, the existence of such a vertex $r$ for a nontrivial block of order at most $k$ implies (i) or (ii), since every nonempty $X_i$ contains at least one vertex.

With this definition in place, we can now state our main structural characterization theorem:

\begin{theorem} For every integer $k \geq 3$, we have $\mathcal{J}_k = \mathcal{H}_k$.
\end{theorem}
\begin{proof} Both classes are closed under disjoint unions, since the parameters are additive over components and the attachment sets are defined within $C_B$. Thus it suffices to consider connected graphs. First, we show that $\mathcal{J}_k \subseteq \mathcal{H}_k$.

Let $G$ be an arbitrary connected graph in $\mathcal{J}_k$. If $G$ is a tree, then it follows immediately that $G \in \mathcal{H}_k$. If $G$ is not a tree, let $B$ denote a nontrivial block of $G$. Since every nontrivial block is 2-connected, we can establish the following general result about 2-connected induced subgraphs of $G$, which will be applied to $B$ itself:

{\noindent\bf Claim 1:} If $N$ is a $2$-connected induced subgraph of $G$, then either $n(N)\le k$ or $n(N)\ge 2k$.

{\noindent \bf Proof of Claim 1:} Assume, for contradiction, that \(k < n(N) < 2k\). Since \(G \in \mathcal{J}_k\), it follows directly from the definition of \(\mathcal{J}_k\) that \(N \in \mathcal{J}_k\). Given \(n(N) < 2k\), we derive \(\nu_k(N) < 2\). However, as \(N\) is 2-connected, the subgraph \(N - x\) remains connected for any vertex \(x \in V(N)\), and \(n(N - x) \ge k\). This implies every generalized \(k\)-vertex cover of \(N\) must contain at least two vertices, so \(\tau_k(N) \ge 2\). But this contradicts the membership of \(N\) in \(\mathcal{J}_k\), which requires \(\nu_k(N) = \tau_k(N)\).
Thus, our initial assumption is false, and we conclude that either \(n(N) \le k\) or \(n(N) \ge 2k\). $\qed$

Applying Claim 1 to the nontrivial block $B$ (which is a 2-connected induced subgraph of $G$), we conclude that either $n(B) \le k$ or $n(B) \ge 2k$. We then have the following claim about the structure of $B$:

{\noindent\bf Claim 2:} If $B$ contains a cycle of length less than $k$, then $n(B) \le k$.

{\noindent \bf Proof of Claim 2:} Suppose, for contradiction, that $B$ contains a cycle of length less than $k$ and $n(B)\ge2k$. By Theorem \ref{ear}, choose an ear decomposition $P_0\cup\cdots\cup P_s$ of $B$ starting with such a cycle. Put $W_i=V(P_0\cup\cdots\cup P_i)$. Each $G[W_i]$ is 2-connected, since it contains the 2-connected ear prefix as a spanning subgraph. Choose the first $j$ with $n(W_j)>k$. Claim 1 gives
\[
n(W_{j-1})\le k,\qquad n(W_j)\ge2k.
\]
Hence the ear $P_j$, with endpoints $u_j,v_j$, has length at least $k+1$.

Choose a shortest $u_j$--$v_j$ path $Q$ in $G[W_{j-1}]$. Its length $q$ satisfies $1\le q\le k-1$, and $C_j=P_j\cup Q$ is a cycle with $\mathrm{dist}_{C_j}(u_j,v_j)=q$. Since $Q$ is induced and $G[W_{j-1}]$ has minimum degree at least two, each endpoint of $Q$ has a neighbor in $W_{j-1}\setminus V(Q)$, hence outside $C_j$. Lemma \ref{cycle-lemma}(iii) implies that $q$ is a multiple of $k$, a contradiction. $\qed$

With the constraint on the size of such blocks established, we now characterize their structural properties relative to the family $\mathcal{H}_k$:

{\noindent\bf Claim 3: } Let $B$ be a nontrivial block of $G$ with $n(B)\le k$. Then
\begin{enumerate}
  \item[(a)] If $n(B) = k$, then $B$ satisfies condition (i) in the definition of $\mathcal{H}_k$.
  \item[(b)] If $n(B) = t < k$, then $B$ satisfies condition (ii) in the definition of $\mathcal{H}_k$.
\end{enumerate}
{\noindent \bf Proof of Claim 3: }(a) Suppose, for contradiction, that $B$ is not an endblock of $G$. Then $B$ contains at least two cut-vertices of $G$, denoted by $u$ and $v$. Since $n(B) = k \ge 3$, there exists a vertex $w \in V(B) \setminus \{u, v\}$. As $B$ is 2-connected, both $B - u$ and $B - v$ are connected. Let $u'$ and $v'$ be neighbors of $u$ and $v$ outside $B$, respectively, and consider the induced subgraph $K = G[V(B) \cup \{u', v'\}]$. Since $G \in \mathcal{J}_k$ and $K$ is an induced subgraph of $G$, it follows that $K \in \mathcal{J}_k$. Note that $n(K) = k + 2 < 2k$ for all $k \ge 3$, so $\nu_k(K) = 1$. However, removing any single vertex from $K$ leaves a connected subgraph of order at least $k$, which implies $\tau_k(K) \ge 2$. This contradicts the definition of $\mathcal{J}_k$, which requires $\nu_k(K) = \tau_k(K)$. Therefore, our initial assumption is false, and $B$ must be an endblock of $G$.

(b) Similarly, suppose, for contradiction, that $B$ contains at least $k-t+2$ cut-vertices of $G$. Let $x_1, x_2, \cdots, x_{k-t+2}$ be cut-vertices of $G$ contained in $B$, and let $y_i$ be a neighbor of $x_i$ outside $B$ for each $i \in \{1, 2, \cdots, k-t+2\}$. Let $K'$ be the induced subgraph $G[V(B) \cup \{y_1, y_2, \cdots, y_{k-t+2}\}]$. Then $n(K') = t + (k-t+2) = k+2$, so $n(K') < 2k$ for all $k \ge 3$. It follows that $\nu_k(K') = 1$. However, removing any single vertex from $K'$ leaves a connected subgraph of order at least $k$, which implies $\tau_k(K') \ge 2$. Since $K'$ is an induced subgraph of $G$ and $G \in \mathcal{J}_k$, we have $K' \in \mathcal{J}_k$, which requires $\nu_k(K') = \tau_k(K')$. This contradiction establishes that $B$ has at most $k-t+1$ cut-vertices in $G$.

Moreover, suppose $B$ contains $m$ cut-vertices where $2\le m\le k-t+1$. Let $a_i=n(X_i)$ and $L=k-t+1$. For contradiction, assume that condition \eqref{small-condition} fails. Then
\[
\sum_{i\ne j}a_i\ge L\qquad(j\in[m]).
\]
Replacing each $a_i$ by $c_i=\min\{a_i,L\}$ preserves these inequalities. Choose integers $0\le b_i\le c_i$ satisfying $\sum_{i\ne j}b_i\ge L$ for every $j\in[m]$, with $s=\sum_i b_i$ as small as possible. If $M=\max_i b_i$, the constraints are equivalent to $s-M\ge L$. If $s-M\ge L+1$, decreasing any positive $b_i$ by one preserves the constraints, a contradiction. Consequently,
\[
s-M=L,\qquad M\le L,\qquad s\le2L.
\]
For each $i$, choose $Y_i\subseteq X_i$ with $n(Y_i)=b_i$ such that $G[Y_i\cup\{x_i\}]$ is connected. This is possible because $G[X_i\cup\{x_i\}]$ is connected. Let
\[
H'=G\left[V(B)\cup\bigcup_{i=1}^mY_i\right].
\]
The graph $H'$ is connected and
\[
n(H')=t+s\le t+2(k-t+1)=2k-t+2\le2k-1,
\]
since $t\ge3$ (see Fig.\ \ref{figure1}). Removing $x_i$ leaves the connected graph formed by $B-x_i$ and all $Y_j$ with $j\ne i$, whose order is at least $t-1+s-b_i\ge t-1+L=k$. Removing a vertex of $B$ other than the $x_i$ leaves a connected graph of at least $t-1+s\ge k$ vertices. Finally, if a vertex of $Y_i$ is removed, then $B$ together with all $Y_j$ for $j\ne i$ remains connected and has at least $t+s-b_i\ge k+1$ vertices. Thus $\tau_k(H')\ge2$, whereas $\nu_k(H')=1$. This contradicts $G\in\mathcal{J}_k$ and proves condition (ii). $\qed$

\begin{figure}[!htbp]
\centering
\begin{tikzpicture}[scale=.85]
\draw[fill=gray!15] (0,0) circle (1.25);
\node at (0,.45) {$B$};\node at (0,.05) {$n(B)=t$};
\coordinate (x1) at (-1.08,-.62);\coordinate (xi) at (0,-1.25);\coordinate (xm) at (1.08,-.62);
\foreach \p in {x1,xi,xm}{\fill (\p) circle (1.8pt);}
\node[left] at (x1) {$x_1$};\node[above left] at (xi) {$x_i$};\node[right] at (xm) {$x_m$};
\draw (x1)--(-3.7,-2) (xi)--(0,-2.5) (xm)--(3.7,-2);
\draw (-3.7,-2.7) ellipse (1.4 and .7);\draw (0,-3.2) ellipse (1.4 and .7);\draw (3.7,-2.7) ellipse (1.4 and .7);
\node at (-3.7,-2.55) {$Y_1$};\node at (-3.7,-2.95) {$n(Y_1)=b_1$};
\node at (0,-3.05) {$Y_i$};\node at (0,-3.45) {$n(Y_i)=b_i$};
\node at (3.7,-2.55) {$Y_m$};\node at (3.7,-2.95) {$n(Y_m)=b_m$};
\node at (-1.85,-2.7) {$\cdots$};\node at (1.85,-2.7) {$\cdots$};
\node at (0,-4.45) {$s=\sum_{i=1}^m b_i,\qquad s-b_i\ge L\quad(i\in[m])$};
\end{tikzpicture}
\caption{Cut-vertices of block $B$ and the selected parts of their attachments.}\label{figure1}
\end{figure}
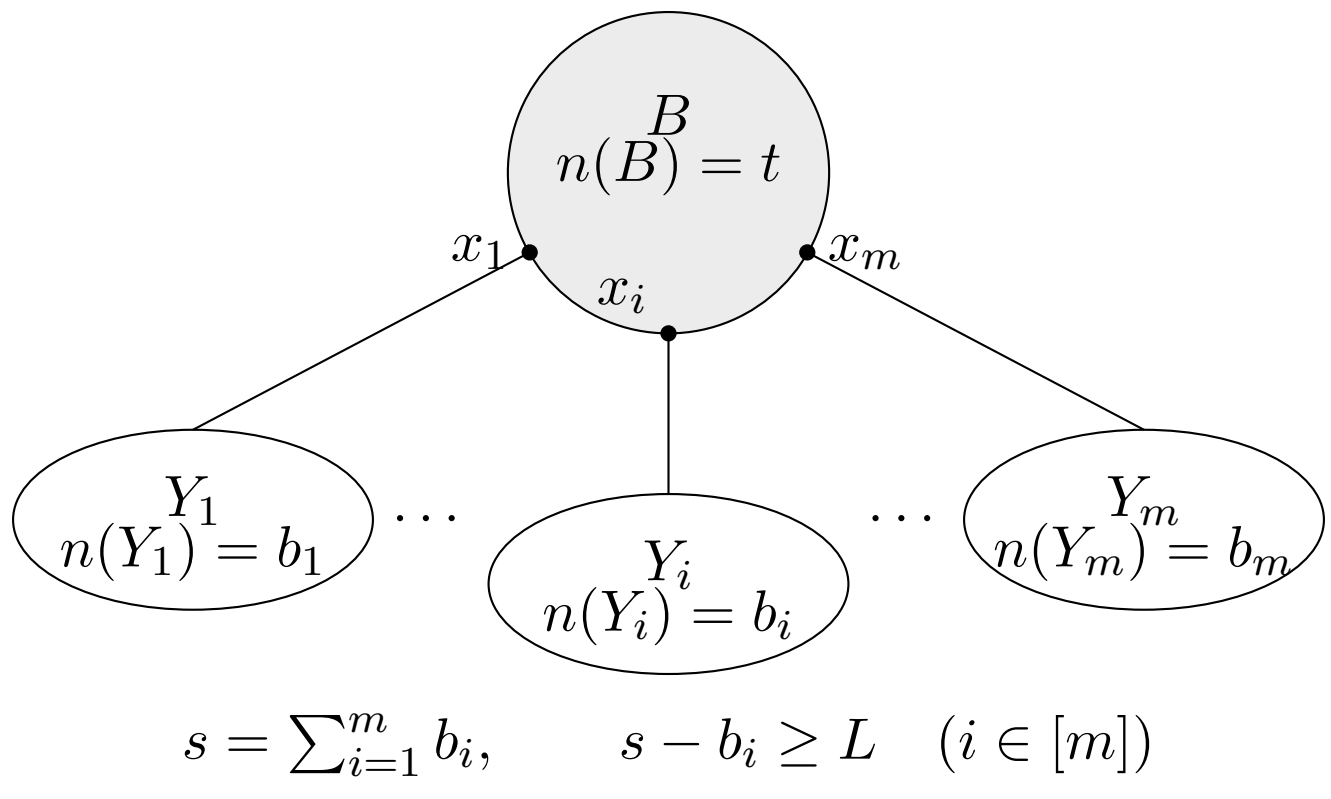

We now consider the remaining case $n(B)>k$. By Claim 2, every cycle in $B$ has length at least $k$, and Lemma \ref{cycle-lemma}(i) implies that every such length is a multiple of $k$. Every cycle in $B$ is chordless. Indeed, a chord would divide a cycle into two cycles whose lengths sum to the original length plus two; all three lengths would be divisible by $k$, implying $k\mid2$, a contradiction.

If $B$ itself is a cycle, Lemma \ref{cycle-lemma}(iii) implies that all its vertices having neighbors outside $B$ occur in one residue class of positions modulo $k$. Choosing every $k$th vertex in this class as an original vertex (or choosing any class if there are no such vertices) gives $B=\text{Sub}_k(H)$ with all cut-vertices of $G$ in $B$ among the vertices of $H$.

Otherwise, split the edges of $B$ into maximal paths $P_0:u_0u_1\cdots u_l$ whose endpoints have degree at least $3$ in $G$ and whose internal vertices have degree $2$ in $G$. Each such path lies on a cycle $C$ of $B$. Since $C$ is chordless, both $u_0$ and $u_l$ have neighbors outside $V(C)$. By Lemma \ref{cycle-lemma}(iii), together with (i), the length $l$ of $P_0$ is a multiple of $k$.

We now define a multigraph $H$ by replacing every such path $u_0u_1\cdots u_{kp}$ by the path $u_0u_k\cdots u_{kp}$ of length $p$. It is clear that $B=\text{Sub}_k(H)$, and every cut-vertex of $G$ contained in $B$ is a vertex of $H$. Thus $G\in\mathcal{H}_k$, and we conclude that $\mathcal{J}_k\subseteq\mathcal{H}_k$.

Next, we establish the inclusion $\mathcal{H}_k \subseteq \mathcal{J}_k$.

First, we verify that $\mathcal{H}_k$ is hereditary. Let $F$ be an induced subgraph of $G\in\mathcal{H}_k$, and let $D$ be a nontrivial block of $F$. The block $D$ is contained in a unique block $B$ of $G$. If $B$ satisfies (i) or (ii), choose $r$ and $U$ as in \eqref{small-region}. If the component of $F$ containing $D$ does not contain $r$, that component lies in $U\setminus\{r\}$ and has fewer than $k$ vertices, so $D$ satisfies (ii). Otherwise, choose $z=r$ if $r\in V(D)$; if $r\notin V(D)$, let $z$ be the unique vertex of $D$ through which the component of $F-V(D)$ containing $r$ attaches to $D$. Uniqueness follows from the maximality of the block $D$. Every vertex outside $U$ can be reached from $D$ only through $r$, hence through $z$. Thus $D$ together with all its attachments in $F$ away from $z$ is contained in $U$ and has order at most $k$. By the equivalence described after \eqref{small-region}, $D$ satisfies (i) or (ii).

If $B=\text{Sub}_k(H)$ satisfies (iii), every internal vertex of a subdivision path has degree $2$ in $G$. If $D$ contains one such vertex, its degree in the 2-connected graph $D$ is at least $2$, so both neighbors and both incident edges lie in $D$. Continuing along the path shows that the entire subdivision path belongs to $D$. Consequently, $D$ is the $k$-subdivision of a connected submultigraph of $H$ containing a cycle. A cut-vertex of $F$ in $D$ has degree at least $3$ in $F$, so it cannot be an internal subdivision vertex. Therefore $D$ satisfies (iii), proving the hereditary property.

It now suffices to show that $\nu_k(G) = \tau_k(G)$ for every connected graph $G \in \mathcal{H}_k$, which we prove by induction on the sum of the order and size of $G$. If $G$ is a tree, Theorem \ref{tree} immediately implies $\nu_k(G) = \tau_k(G)$. If $G$ is a cycle of length less than $k$, both parameters are zero. If its length is a multiple of $k$, then $\nu_k(G)=\tau_k(G)$ holds as well. These are the only cycle lengths permitted by the definition of $\mathcal{H}_k$. We may therefore assume that $G$ is neither a tree nor a cycle.

Suppose that $G$ contains a block $B$ satisfying condition (i). Then $\nu_k(G) = \nu_k(G - V(B)) + 1$ and $\tau_k(G) = \tau_k(G - V(B)) + 1$. By the hereditary property established above, $G - V(B) \in \mathcal{H}_k$. Applying the induction hypothesis, we have $\nu_k(G - V(B)) = \tau_k(G - V(B))$. Combining these equalities yields $\nu_k(G) = \tau_k(G)$.

Suppose $G$ contains a block $B$ satisfying condition (ii). Choose $r$ and $U$ as in \eqref{small-region}, and let $G'$ be obtained by replacing $B$ by a spanning tree on the same vertex set. The graph $G'[U]$ remains connected, and $G'$ has fewer edges than $G$. We show that $\nu_k(G')=\nu_k(G)$ and $\tau_k(G')=\tau_k(G)$.

In both graphs, every $k$-tree meeting $U$ contains $r$, since $n(U\setminus\{r\})<k$ and this set has no neighbors outside $U$. Hence at most one member of a generalized $k$-matching meets $U$. If a $k$-tree $T$ in such a matching meets $U$, remove its vertices in $U\setminus\{r\}$. The remaining tree is connected and contains $r$: no path in $T$ can enter and leave $U\setminus\{r\}$ without passing through $r$ twice. In the other graph, grow a tree from $r$ inside $U$ using the same number of additional vertices, and join it to the remaining part of $T$. This gives a $k$-tree disjoint from all other members of the matching. Applying this argument in both directions proves $\nu_k(G')=\nu_k(G)$.

If a generalized $k$-vertex cover $S$ of $G$ meets $U$, then $(S\setminus U)\cup\{r\}$ covers $G'$ and has cardinality at most $|S|$. If $S$ is disjoint from $U$, replacing the connected graph on $U$ by another connected graph on the same vertex set leaves the vertex sets of the components after deleting $S$ unchanged. Thus $S$ also covers $G'$. The same argument with $G$ and $G'$ interchanged proves $\tau_k(G')=\tau_k(G)$.

The nontrivial blocks of $G'$ are precisely the nontrivial blocks of $G$ other than $B$. Their attachment sets are unchanged, since the replacement preserves connectivity between all vertices of $B$; at a cut-vertex shared with another block, the vertices on the $B$ side remain on that side. Thus $G'\in\mathcal{H}_k$. By induction, $\nu_k(G')=\tau_k(G')$, and consequently $\nu_k(G)=\tau_k(G)$.

In all such cases, the induction hypothesis yields $\nu_k(G) = \tau_k(G)$. We may therefore assume that $G$ contains no such blocks. Let $B = \text{Sub}_k(H)$ for some multigraph $H$, where every cut-vertex of $G$ contained in $B$ is a vertex of $H$. The multigraph $H$ is connected and contains a cycle. Every internal subdivision vertex has degree $2$ in $G$, so each component of $G-V(H)$ meeting $V(B)\setminus V(H)$ is an internal subdivision path of order $k-1$. Thus the additional hypothesis of Lemma \ref{multigraph} holds. Applying that lemma to $B$, we obtain $\nu_k(G) = \nu_k(G - V(H)) + n(H)$ and $\tau_k(G) = \tau_k(G - V(H)) + n(H)$. Since $\mathcal{H}_k$ is hereditary, the induction hypothesis gives $\nu_k(G - V(H)) = \tau_k(G - V(H))$, which implies $\nu_k(G)=\tau_k(G)$. Since $\mathcal{H}_k$ is hereditary, the same equality holds for every induced subgraph of $G$, proving $G\in\mathcal{J}_k$ and completing the proof. $\qed$
\end{proof}





\end{spacing}

\end{document}